\documentclass[11pt]{amsart}
\usepackage{txfonts}
\usepackage{}
\usepackage{amssymb}
\theoremstyle{plain}
\newtheorem{Thm}{Theorem}

\begin{document}

\title[K\"ahler-Ricci flow on complete non-compact manifolds]
{Global K\"ahler-Ricci flow on complete non-compact manifolds }

\author{Ma Li}
\address{Department of mathematics \\
Henan Normal university \\
Xinxiang, 453007 \\
China} \email{lma@tsinghua.edu.cn}

\thanks{The research is partially supported by the National Natural Science
Foundation of China (N0.11271111)}

\begin{abstract}
In this paper, we study the global K\"ahler-Ricci flow on a complete
non-compact K\"ahler manifold. We prove the following result. Assume
that $(M,g_0)$ is a complete non-compact K\"ahler manifold such that
there is a potential function $f$ of the Ricci tensor, i.e.,
$$
R_{i\bar{j}}(g_0)=f_{i\bar{j}}.
$$
Assume that the quantity $|f|_{C^0}+|\nabla_{g_0}f|_{C^0}$ is finite
and the $L^2$ Sobolev inequality holds true on $(M,g_0)$. Then the
Kahler-Ricci flow with the initial metric $g_0$ either blows up at
finite time or infinite time to Ricci flat metric or exists globally
with Ricci-flat limit at infinite time. Related results are also
discussed.


\end{abstract}

 \maketitle

\section{Introduction}\label{sect1}
In the study of Ricci flow, the main problem is to understand the
global behavior of it \cite{H2}. Since one may meet singularities of
Ricci flow in finite time, it is interesting to know the singularity
structure of the Ricci flow. R.Hamilton conjectured that after
suitable normalization, the Ricci flow at the blow-up point looks
like a Ricci soliton. In this paper, we study the global behavior
of the K\"ahler-Ricci flow on a complete non-compact manifold with
bounded curvature and we show in this case Hamilton's conjecture is
right. We remark that the local existence of the flow has been
proved by W.Shi \cite{sh} and we use this local flow in this paper.

We prove the following result.

\begin{Thm}\label{main}
Assume that $(M^n,g_0)$ ($n=dim_CM$) is a complete non-compact
K\"ahler manifold with bounded curvature and there is a potential
function $f$ of the Ricci tensor, i.e.,
$$
R_{i\bar{j}}(g_0)=f_{i\bar{j}}.
$$
Assume the quantity $ |f|_{C^0}+|\nabla_{g_0}f|_{C^0}$ is finite. We
also assume that the $L^2$ Sobolev inequality holds true on
$(M,g_0)$. Then the Kahler-Ricci flow
\begin{equation}\label{KF}
\partial_tg_{i\bar{j}}=-R_{i\bar{j}}(g)
\end{equation}
with the initial metric $g_0$ either blows up at finite time or
infinite time to Ricci flat metric or exists globally with
Ricci-flat limit at infinite time.
\end{Thm}

Our assumptions are natural. Note that in the statement of our result, we
assume $L^2$ Sobolev inequality, namely, there exists a global
constant $C>0$ such that for any $u\in C_0^1(M,g_0)$,
$$
C(\int_M |u|^{2n/(n-1}dV)^{(n-1)/n}\leq
\int_M(|\nabla_{g_0}u|^2+u^2)dV),
$$
 which, plus the bounded curvature assumption, implies the Gross Log-Sobolev inequality and the non-local collapsing for the metric $g_0$.
 Hence the W-functional introduced by Perelman \cite{P} makes sense on the Riemannian manifold $(M,g_0)$.
 Along the Ricci flow, as long as the curvature is bounded, all the Riemannian metrics are equivalent each other before the first singularity.
 This implies that the W-functional is well-defined along the Ricci flow before it blows up. In our proof of
Theorem \ref{main}, we shall use the the equivalent inequality for
the W-functional (and see \cite{CN} for more discussions).

Without the non-local collapsing for the metric $g_0$, we have the
following result.

\begin{Thm}\label{main+1}
Assume that $(M,g_0)$ is a complete non-compact K\"ahler manifold
with bounded curvature and there is a potential function $f$ of the
Ricci tensor, i.e.,
$$
R_{i\bar{j}}(g_0)=f_{i\bar{j}}.
$$
Assume only $ |\nabla_{g_0}f|_{C^0}$ is finite and the initial
metric has non-negative bisectional curvature. There is a global
Kahler-Ricci flow with the initial metric $g_0$.
\end{Thm}

To understand the global behavior of the K\"ahler-Ricci flow
established by W.Shi \cite{sh}, we need to use some invariants for
the various curvature quantities. The key problem is weather the positivity of the Ricci
curvature under the K\"ahler-Ricci flow holds true, which is central for understanding
the convergence of the K\"ahler-Ricci flow on K\"ahler-Einstein
manifolds. It is well-known that the K\"ahler-Ricci flow preserves
the positivity of the curvature operator, the positivity of the
bisectional curvature and the positivity of the scalar curvature.
 Phong and Sturm
\cite{phong} have found that the positivity of the Ricci curvature
is preserved by the K\"ahler-Ricci flow on compact K\"ahler
manifolds of dimension two, under the assumption that the sum of any
two eigenvalues of the traceless curvature operator on traceless
(1,1)-forms is non-negative. Since their argument only uses the
maximum principle, one may easily extend their result to a complete
noncompact manifold with the same curvature assumptions above. As
applications of our Theorem \ref{main}, we can prove that in either
cases as above, there is a global K\"ahler-ricci flow on the K\"ahler
manifold with Ricci flat metric as its limit at time infinity.

We remark that one may formulate more convergence results using the
invariant set constructed by B.Wilking \cite{W}, who finds almost
all invariant curvature conditions of Ricci flows. See also related interesting works \cite{SW} \cite{H}
\cite{CT} \cite{Ch}\cite{C}\cite{BS} \cite{N} \cite{Ni} \cite{CT2} \cite{Z}.

The plan of the paper is below. In section \ref{sect2}, we consider
the non-local collapsing result obtained by using Perelman's
W-functional and the scalar
curvature bound found by B.Chow in our K\"ahler-Ricci
flow case. In sections \ref{sect3} and \ref{sect4}, we prove Theorem
\ref{main}. We also consider some consequences of our Theorem \ref{main}
in the last section.

\section{non-local collapsing and Perelman's
W-functional and uniform scalar curvature bound}\label{sect2}

We now recall some facts about Perelman's W-functional on a closed
Riemannian manifold $(M^n,g)$. First G.Perelman \cite{P} defines the
F-functional as
$$
F(g,f)=\int_M (R+|\nabla f|^2)dm
$$
where $dm=e^{-f}dv_g$ is a fixed measure and $R$ is the scalar
curvature of the metric $g$. This functional is an extension of
Hilbert action on the Riemannian metric space on $M$. Let $v$ and
$h$ be the variations of $g$ and $f$ respectively. Since $dm$ is
fixed, we have the relation $h=tr_hv/2$.
Then we have
$$
\delta F(v,h)=-\int(v,Rc(g)+D^2f)dm
$$
where $Rc(g)$ is the Ricci tensor of $g$ and $D^2f$ is the Hessian
matrix of the function $f$. For any fixed parameter $\tau>0$, the W-functional is defined by
$$
W(g,f)=W(g,f,\tau)=\int_M[\tau(R+|\nabla
f|^2)+f-n](4\pi\tau)^{-n/2}dm.
$$
That is,
$$
W(g,f)=(4\pi\tau)^{-n/2}\tau F(g,f)+(4\pi\tau)^{-n/2}\int_M(f-n)dm.
$$
The importance of W-functional is that it is a diffeomorphism
invariant. Then we have
$$
\delta
W(v,h)=-\tau\int_M(v,Rc(g)+D^2f-\frac{1}{2\tau}g)(4\pi\tau)^{-n/2}dm.
$$
Hence, the $L^2$ gradient flow of $W(g,f)$ is
$$
\partial_tg=-2(Rc(g)+D^2f-\frac{1}{2\tau}g)
$$
with
$$
f_t=-\Delta_gf-R+\frac{n}{2\tau}.
$$

Let $\phi(t)$ be the flow generated by the time-dependent
vector-field $\nabla_gf$. Let $\bar{g}(t)=\phi(t)^*g(t)$ and
$\bar{f}(t)=\phi(t)^*f(t)$. Then we have
$$
\partial_t\bar{g}=-2(Rc(\bar{g})-\frac{1}{2\tau}\bar{g})
$$
with
$$
\bar{f}_t=-\Delta_{\bar{g}}\bar{f}-\bar{R}+|\nabla_{\bar{g}}\bar{f}|^2+\frac{n}{2\tau}.
$$
Let $C(s)=1-\frac{s}{\tau}$,
$$
t(s)=-\tau\log C(s).
$$
and
$$
\check{g}(s)=C(s)\bar{g}(t(s)).
$$
Then we have
$$
\partial_t\check{g}=-2Rc(\check{g})
$$
and
$$
\check{f}_t=-\Delta_{\check{g}}\check{f}-\check{R}+|\nabla_{\check{g}}\check{f}|^2+\frac{n}{2\tau
C(s)}.
$$
Define
$$
\mu(g,\tau)=\inf \{W(g,f,\tau), \int_Mdm=1\}.
$$
Then for any fixed $\tau>0$, along the Ricci flow, $\mu(g(t),\tau)$
is non-decreasing.

We assume that the scalar curvature of the K\"ahler-Ricci flow is
uniformly bounded, which can be proved in our situation (see section
\ref{sect3}) as in Theorem \ref{main}. In a K\"ahler manifold
$(M^n,g)$ (where $n$ is the complex dimension), we normalize the
W-functional as
$$
W(g,f,\tau)=(4\pi\tau)^{-n}\int_M[2\tau(R+|\nabla
f|^2)+f-2n]e^{-f}dV
$$
with $(4\pi\tau)^{-n}\int_Me^{-f}dV=1$ and define
$$
\mu(g,\tau)=\inf_{f|(4\pi\tau)^{-n}\int_Me^{-f}dV=1} W(g,f,\tau)
$$
Assume now that $(M,g)$ is a complete noncompact Riemannian
manifold.  Note that for $u=e^{-f/2}$, we have
$$ (4\pi\tau)^{-n}\int_M u^2dV=1
$$
and
$$
W(g,f,\tau)=(4\pi\tau)^{-n}\int_M[2\tau(Ru^2+4|\nabla u|^2)+u^2\log
u^2-2nu^2]dV
$$
 which will be written as $W(g,u,\tau)$. This latter formulation
 allows us to define the W-functional for $u\in C_0^\infty(M)$ (that is $e^{-f/2}\in C_0^\infty(M)$ and this is noticed
by G.Perelman, see Remark 3.2 in \cite{P}, see also p.240 of
\cite{CN}) and
 this fact is very useful in analyzing the properties of
 $W$-functional. Applying the decay estimates of Hamilton \cite{H2} and Dai-Ma \cite{DM},
 Li-Yau-Hamilton gradient estimate \cite{CN}, and Perelman's Proposition 9.1 \cite{P},
 we can show that W-functional is also non-decreasing along the
 Ricci flow on a complete non-compact Riemannian manifold with
 the initial bounded curvature. Since the argument is well-known to
 experts, we omit the detail.

The important application about W-functional is the no local
collapsing result due to G.Perelman: \emph{Using this property and
the scalar curvature bound of the Ricci flow on $(0,T)$ with
$T<\infty$, we can get the no local-collapsing at $T$}, and the no
collapsing data depends only on the initial metric and the global
scalar curvature bound. For the full proof of this result, one may
refer to G.Perelman's paper \cite{P}. See also Sesum-Tian's paper
\cite{ST} on compact Kahler manifolds. We remark that the argument
can be carried out to our complete non-compact case. Since the
argument is almost the same, we omit the detail.

We now recall some formulae from B.Chow's paper \cite{chow}. Along the K\"ahler-Ricci flow (\ref{KF}), we have
$$
R_{i\bar{j}}=-\partial_i\partial_{\bar{j}}\log
det(g_{k\bar{l}})=\partial_i\partial_{\bar{j}}f,\; \; R=\Delta f,
$$
and
$$
\partial_tR_{i\bar{j}}=-
\nabla_i\nabla_{\bar{j}}(g^{k\bar{l}}\partial_tg_{k\bar{l}})=\nabla_i\nabla_{\bar{j}}R,
$$
which implies that
$
\partial_i\partial_{\bar{j}}(\partial_t-\Delta)f=0.
$
We can normalize $f$ such that it satisfies $
\partial_t f=\Delta f $
with the uniform estimate $t|\nabla f|^2+f^2\leq C$,
where $C$ is some uniform constant depending only on the initial
data $f_0$ and $|\nabla f|^2=g^{k\bar{l}}f_kf_{\bar{l}}$.

Recall that the scalar curvature evolution equation is
\begin{equation}\label{scalar}
(\partial_t-\Delta)R=|R_{i\bar{j}}|^2\geq \frac{1}{n}R^2.
\end{equation}
Applying the maximum principle we get the lower bound $R\geq
-\frac{n}{t}$.
We remark that we can use the maximum principle for the
quantities below before the first singularity of the Ricci flow
since we always use the local solution obtained by Shi \cite{sh}
with bounded curvature. The following result is from B.Chow's paper \cite{chow}.
\begin{Thm}\label{bound-1} Assume that there is a potential function $f$ on $M$
such that
$$
R_{i\bar{j}}(g_0)=f_{i\bar{j}}
$$
with bounded gradient $|\nabla_{g_0}f|_{C^0}$. Then there is a
smooth function $f=f(t)$ such that $$
R_{i\bar{j}}(g(t))=f(t)_{i\bar{j}}
$$
and we also have some uniform constant $C_0>0$ depending only on
the initial metric such that
$$
R+|\nabla_{g(t)} f|^2\leq C_0.
$$
\end{Thm}

The proof of Theorem \ref{bound-2} below is similar to Chow's argument \cite{chow},
but some new barrier function should be constructed.

\begin{Thm}\label{bound-2} Assume the same condition for the global Kahler-Ricci flow
$g(t)$ with the potential function $f$ being a bound function. Then
we also have some uniform constant $C_1>0$ depending only on the
initial metric and the dimension $n$ such that
$$
R+|\nabla_{g(t)} f|^2\leq C_1/t.
$$
\end{Thm}

\begin{proof} (of Theorem \ref{bound-2})
By direct computation \cite{chow}, we have
\begin{equation}\label{star+1}
(\partial_t-\Delta) [t(|\nabla f|^2+R)]\leq
-\frac{2t}{n}R^2+R+2|\nabla f|^2.
\end{equation}
Note that
$$
\frac{t}{n}R^2=\frac{t}{n}(R+2|\nabla f|^2)^2-\frac{4t|\nabla
f|^2}{n}(R+|\nabla f|^2).
$$
Then we get from (\ref{star+1}) that
\begin{equation}\label{star+2}
 \begin{array}{l l}
(\partial_t-\Delta) [t(|\nabla f|^2+R)] &\leq
-\frac{t}{n}(R+2|\nabla f|^2)^2 \\
&-\frac{t}{n}R^2+R+2|\nabla f|^2+\frac{4t|\nabla f|^2}{n}(R+|\nabla
f|^2).
 \end{array}
\end{equation}
Assume that $R\geq 0$, using $|\nabla f|^2\leq \frac{C}{t}$, we
obtain from (\ref{star+2}) that
\begin{equation}\label{star+3}
 \begin{array}{l l}
(\partial_t-\Delta) [t(|\nabla f|^2+R)] &\leq
-\frac{t}{n}(R+2|\nabla f|^2)^2 \\
&-\frac{t}{n}R^2+C(R+2|\nabla f|^2).
 \end{array}
\end{equation}
Note that the term $R+2|\nabla f|^2$ is bounded by $C_0$, and if
$R\geq 0$ and  $\frac{t}{n}(R+2|\nabla f|^2)\geq 2CC_0$, we then
have
$$
(\partial_t-\Delta) [t(|\nabla f|^2+R)] \leq 0.
$$
As in \cite{sh}, any $T>0$, we construct a distance like positive
function $h(x)$ such that for any $t\in [0,T]$, $
\Delta_{g(t)} h\leq 1$,
and
$
h(x)\to \infty
$
as $d(x,\bar{x})\to\infty$, where $\bar{x}$ is a fixed point in $M$.
Define the domain
$$\Omega_T=\{(x,t)\in M\times (0,T], R\geq 0, \frac{t}{n}(R+2|\nabla f|^2)\geq 2CC_0 \}.
$$
Note that in $\Omega_T$, for $\epsilon>0$,
$$
\begin{array}{l l}
&(\partial_t-\Delta) [t(|\nabla f|^2+R)+2\epsilon t-\epsilon
h(x)]\\
&= (\partial_t-\Delta) [t(|\nabla f|^2+R)]-2\epsilon+\epsilon
\Delta_{g(t)}h(x)\leq -\epsilon<0.
\end{array}
$$
Since the region $U_\epsilon$ for the function $t(|\nabla
f|^2+R)-2\epsilon t-\epsilon h(x)\geq 2CC_0$ is compact and using
the maximum principle in $U_\epsilon$, we know that
$$
t(|\nabla f|^2+R)\leq 2CC_0+2\epsilon t+\epsilon h(x), \ \ (x,t)\in
U_\epsilon \cap \Omega_T.
$$
Note that $\Omega_T=\bigcup_{\epsilon>0} U_\epsilon$ and $U_\epsilon
\subset U_{\epsilon'}$ for $\epsilon>\epsilon'$. Then for any
$(x_0,t_0)\in \Omega_T$, we have some $\epsilon_0>0$ such that
$(x_0,t_0)\in U_\epsilon $
for any $0<\epsilon<\epsilon_0$ and $$t(|\nabla f|^2+R)(x_0,t_0)\leq
2CC_0+ 2\epsilon t_0+\epsilon h(x_0).$$ Letting $\epsilon\to 0$, we
have
$$
t(|\nabla f|^2+R)(x_0,t_0)\leq 2CC_0.
$$
That is to say, in $\Omega_T$, we have
$ t(|\nabla f|^2+R)\leq 2CC_0. $
Therefore, for any $T>0$, we always have $|R|\leq 2\max(CC_0,n)/t$
in $M\times (0,T)$. This completes the proof of Theorem \ref{bound-2}.

\end{proof}

\section{global flow and finite time blow-up is Ricci flat}\label{sect3}

In this section, we shall prove Theorem \ref{main+1}.
Using our recent result about extension result of Ricci flow
\cite{MC}, which extends the beautiful result of N.Sesum on compact
manifold, we know that if the K\"ahler-Ricci flow blows up at first
finite time $T$, then the Ricci curvature must blow up at time $T$.
Using the scalar curvature bound of the K\"ahler-Ricci flow in
previous section, we can show the Ricci-flatness of the blow up
limit of the flow.

\begin{Thm}\label{finite} Assume the conditions as in Theorem
\ref{main}. If the Kahler-Ricci flow exists only in finite time,
then its blow-up limit at the maximal time exists and the limit is
Ricci flat.
\end{Thm}

We remark that the the blow-up limit or just limit in the results
above and below is only for some subsequence, not for full blow-up
sequence. The process to prove Theorem \ref{finite} is below. First we note
that by using Perelman's non-local-collapsing result at any finite
time, we can applying Hamilton's Cheeger-Gromov convergence theorem \cite{CN}
to take the blow-up limit $g_\infty(t)$. Then we use the uniform
scalar curvature bound of the original Ricci flow to obtain
$R(g_\infty)=0$, which implies the Ricci flatness of the limit
metric $g_\infty$.

\begin{proof} (for Theorem \ref{finite}).
Let $T<\infty$ be the first blow-up time of the Ricci flow
$(M,g(t))$. Choose some constants $C>0$, the finite time $t_j\to T-$
and take points $x_j$ such that
$$
\sup_{M\times [t_j-C^{-1}K_j^{-1},t_j]}|Rm(g(t))|\leq CK_j
$$
where $K_j=|Rm(g(t_j))|(x_j)\to\infty$. Define
$$
g_j(t):=K_jg(t_j+K_j^{-1}t).
$$
Clearly we have $|Rm(g_j(t))|\leq C$ and $|R(g_j(t))|\leq
C/K^{-1}_j\to 0$ . Using Perelman's non-local collapsing result we
know the injectivity radius bound condition along the flow (\ref{KF}). Then we can get the
limit Ricci flow $g_\infty(t)$ of $g_j(t)$ with
$|Rm(g_\infty(t))|\leq C$ and $R(g_\infty(t))=0$ for $t\in
(-\infty,\omega)$. Using the scalar curvature evolution equation (\ref{scalar})
we get that $Rc(g_\infty(t))=0$. This completes the proof of Theorem
\ref{finite}.
\end{proof}

We now prove Theorem \ref{main+1}: Assume that we have a local
K\"ahler-Ricci flow with initial metric $g_0$ with bounded curvature
and non-negative bisectional curvature. Using the result of
Bando-Mok, we know that the non-negativity bisectional curvature is
preserved along the flow. Hence the Ricci curvature is non-negative
along the Ricci flow. Using the scalar curvature bound and
non-negativity Ricci curvature condition, we get via the use of the
elementary relation $R_{i\bar{j}}\leq Rg_{i\bar{j}}$ the uniform
bound for the Ricci curvature. Then using Theorem 1.4 \cite{MC}, we
know that the K\"ahler-Ricci flow exists globally with bounded
curvature.

\section{global Kahler-Ricci flow and its limit at $t=\infty$}\label{sect4}

In this section, we prove Theorem \ref{main}. Note that there is no
local collapsing for the Ricci flow at any finite time according to
the result of perelman mentioned in section two. If the curvature of
the flow blows up at infinite time, we then can take the similar
blow-up limit as in the finite time case \cite{MC} and show that its
blow-up limit is Ricci flat. Otherwise, the curvature of the flow is
uniformly bounded and by Theorem 4, we can simply take the limit,
which is scalar flat at first. Using the evolution equation (\ref{scalar})
again we know that the limit is also Ricci flat.
This completes the proof of Theorem \ref{main}.

We give two applications of Theorem \ref{main}.
  One is for the metric $g_0$ with non-negative bisectional
  curvature, which implies the non-negativity of the Ricci
  curvature. Note that the non-negativity of the bisectional
  curvature property is preserved by Kahler-Ricci flow as proved by Bando-Mok as mentioned in the introduction.
  Note that in this case, the scalar curvature dominates other curvatures.
  Since we have global bound of the scalar curvature of the Ricci flow, we must have the global K\"ahler-Ricci flow
  with bounded curvature. Hence by Theorem \ref{main}, the flow is
  global one and has its limit at $t=\infty$ the Ricci flat metric.
Then we have the following result, which may also be derived from
the arguments or results in \cite{NT} and \cite{PSW}.
\begin{Thm}Assume that $(M,g)$ is a complete non-compact K\"ahler manifold with bounded curvature
such that there is a potential function $f$ of the Ricci tensor,
i.e.,
$$
R_{i\bar{j}}(g_0)=f_{i\bar{j}}.
$$
Assume the quantity $ |f|_{C^0}+|\nabla_{g_0}f|_{C^0}$ is finite.
Suppose the initial metric has non-local-collapsing and the $L^2$
Sobolev inequality holds true on $(M,g_0)$. Assume also that the
initial metric has non-negative bisectional curvature. Then the
K\"ahler-Ricci flow with the initial metric $g_0$ exists globally with
Ricci-flat limit at infinite time.
\end{Thm}

  The other is for complex dimension two, where we assume that
  $(M,g)$ is a complete non-compact K\"ahler surface with the bounded curvature, positive Ricci curvature, and
  positive partial curvature operator $S$ in the case handled by Phong-Sturm
  \cite{phong}. We make this precise and define the partial curvature
  operator
$$
S_{i\bar{j}k\bar{l}}=R_{i\bar{j}k\bar{l}}-\frac{1}{n}(S_{i\bar{j}}\delta_{k\bar{l}}
+\delta_{i\bar{j}}S_{k\bar{l}})+\frac{1}{n^2}R\delta_{i\bar{j}}\delta_{k\bar{l}}
$$
acting on the space of traceless $(1,1)$ forms on $M$. Here,
$R_{i\bar{j}k\bar{l}}$ is the curvature tensor of $g$,
$S_{i\bar{j}}=R_{i\bar{j}}-\frac{1}{n}Rg_{i\bar{j}}$.

Then we have
\begin{Thm} Assume that $(M,g)$ is a complete non-compact K\"ahler surface with bounded curvature
such that there is a potential function $f$ of the Ricci tensor,
i.e.,
$$
R_{i\bar{j}}(g_0)=f_{i\bar{j}}.
$$
Suppose the quantity $ |f|_{C^0}+|\nabla_{g_0}f|_{C^0}$ is finite.
Assume that the initial metric has non-local-collapsing and the
$L^2$ Sobolev inequality holds true on $(M,g_0)$. Assume also that
the initial metric has non-negative non-trivial Ricci curvature and
the non-negative sum of the two lowest eigenvalues of the operator
$S$ acting on the space of traceless $(1,1)$ forms on $M$. Then the
K\"ahler-Ricci flow with the initial metric $g_0$ exists globally with
Ricci-flat limit at infinite time.
\end{Thm}

{\bf Acknowledgement}. The author thanks the unknown referees very much for useful suggestions.

\end{document}